\documentclass[12pt]{article}

\usepackage{amsfonts}
\usepackage{amssymb}
\usepackage{amsthm}
\usepackage{amsmath}
\usepackage{newlfont}
\usepackage{graphicx}
\usepackage[all,arc]{xy}
\usepackage{epsfig}
\usepackage{amsmath,amscd}
\usepackage{url}

\usepackage{multicol}
\usepackage[T1]{fontenc}

\newtheorem {proposition}{Proposition}[section]
\newtheorem {theorem}{Theorem}[section]
\newtheorem {lemma}{Lemma}[section]
\newtheorem {example}{Example}[section]
\newtheorem {definition}{Definition}[section]

\newtheorem {corollary}{Corollary}[section]

\newtheorem{remark}{Remark}

\author{{\DJ{}or\dj{}e Barali\' {c}}\\ {\small Mathematical Institute SASA}\\[-2mm] {\small Belgrade, Serbia}
 \and Rade \v{Z}ivaljevi\'{c}\\ {\small Mathematical Institute SASA}\\[-2mm] {\small Belgrade, Serbia}}

\title{Colorful versions of the\\ Lebesgue, KKM, and Hex theorem}
\date{}
\begin{document}
\maketitle

\begin{abstract}\noindent
Following and developing ideas of  R. Karasev (Covering dimension
using toric varieties, arXiv:1307.3437), we extend the Lebesgue
theorem (on covers of cubes) and the
Knaster-Kuratowski-Mazurkiewicz theorem (on covers of simplices)
to different classes of convex polytopes (colored in the sense of
M.~Joswig).  We also show that the $n$-dimensional Hex theorem
admits a generalization where the $n$-dimensional cube is replaced
by a $n$-colorable simple polytope. The use of specially designed
{\em quasitoric manifolds}, with easily computable cohomology
rings and the cohomological cup-length, offers a great flexibility
and versatility in applying the general method.
\end{abstract}

\renewcommand{\thefootnote}{}
\footnotetext{The authors were supported by the Grant 174020 of
the Ministry for Education and Science of the Republic of Serbia.
}

\section{Introduction}

The well known connection between the classical
Lyusternik--Schnirelmann category (LS-category) and the
cohomological cup-length is a simple, yet elegant and powerful
method of studying geometric/topological properties of a space by
computable invariants arising in algebraic topology. Together with
its generalizations and ramifications, this connection is indeed
one of evergreen themes of geometry and topology.

\medskip
It was an interesting recent observation of Karasev \cite{Karasev}
that a similar cohomological cup-length approach can be utilized
for the proof of some results of more combinatorial nature,
including the following two classical results of Lebesgue, and
Knaster, Kuratowski, Mazurkiewicz (KKM).

\begin{theorem} \label{thm:Lebesgue}
{\em (Lebesgue)} If the unit cube $[0, 1]^n$ is covered by a
finite family  $\{X_i\}_{i\in I}$ of closed sets so that no point
is included in more than $n$ sets, then one of them must intersect
two opposite facets of the cube.
\end{theorem}

\begin{theorem}\label{thm:KKM}
{\em (KKM)} If a non-degenerate simplex $\Delta^n\subset
\mathbb{R}^n$ is covered by a finite family $\{F_i\}_{i\in I}$ of
closed sets so that no point is covered more than $n$ times then
one of the sets $F_i$ intersects all the facets of $\Delta^n$.
\end{theorem}

The method of Karasev was based on the use of cohomological
properties of (both non-singular and singular) toric varieties. In
particular he was able to unify Theorems~\ref{thm:Lebesgue} and
\ref{thm:KKM} and interpret them as special cases of a single
statement valid for all simple polytopes.

\begin{theorem} {\em \cite[Theorem~5.2.]{Karasev}}
Suppose that a simple polytope $P\subset \mathbb{R}^n$ is covered
by a family of closed sets $\{X_i\}_{i\in I}$ with covering
multiplicity at most $n$. Then for some $i\in I$ the set $X_i$
intersects at least $n + 1$ facets of $P$.
\end{theorem}

We continue this study by methods of toric topology, emphasizing
the role of quasitoric manifolds and Davis-Januszkiewicz spaces
\cite{Davis, BuPan}. We focus on special classes of simple
polytopes including the class of $n$-colorable simple polytopes
which were introduced by Joswig in \cite{Joswig}. The associated
classes of quasitoric manifolds have computable and often
favorable cohomological properties, which have already found
applications outside toric topology  \cite{Bar, BarGru}.

\medskip
Our central results, the `Colorful Lebesgue theorem'
(Theorem~\ref{thm:color-Lebesgue}) and the `Colorful KKM-theorem
(Theorem~\ref{thm:color-KKM}), together with their companions
Theorem~\ref{thm:quantitative-Lebesgue} and
Theorem~\ref{thm:color-KKM-bis}, are designed to include
Theorems~\ref{thm:Lebesgue} and \ref{thm:KKM} as special cases and
to illuminate the role of special classes of quasitoric manifolds
over $n$-colorable and $(n+1)$-colorable simple polytopes.

\medskip
In the same vein we prove the `Colorful Hex theorem'
Theorem~\ref{thm:color-Hex} and describe a `Colorful Voronoi-Hex
game' played by $n$ players on an $n$-dimensional Voronoi
checkerboard.

\medskip
We referrer the reader, curious or intrigued by the use of the
word `colorful' in these statements, to \cite{Aro} and
\cite{Ziegler} for a sample of results illustrating how the term
`colorful' gradually acquired (almost) a technical meaning in many
areas of geometric and topological combinatorics.

\section{Overview and preliminaries}

A basic insight from the theory of  {\em Lebesgue covering
dimension} is that an $n$-dimensional space cannot be covered by a
family $\mathcal{U}$ of open sets which are {\em `small in size'}
unless we allow non-empty intersections of $(n+1)$ sets or more.
In other words if know in advance that the covering multiplicity
of the family is $\leq n$, then some of the sets $U\in
\mathcal{U}$ must be `large' in some sense.

\medskip
Theorems~\ref{thm:Lebesgue} and \ref{thm:KKM} turn this vague
sense of `largeness' into precise results where the combinatorics
and facial structure of the cube and simplex respectively plays an
important role.

\medskip
Karasev \cite{Karasev} has found a very natural and interesting
way of proving and generating such results, based on the theory of
(complex and real) toric varieties. The use of the cup-length
estimates is of course well known in the theory of
Lyusternik-Schnirelmann category and its ramifications. However,
more combinatorial aspects of the problem and possibilities of the
method don't seem to have been carefully explored and they
certainly deserve a further study.

\medskip
From this point of view it is quite natural to explore which
classes of convex polytopes may provide an adequate concept of
`largeness' suitable for generalizing classical theorems of
Lebesgue (on coverings of cubes) and the
Knaster-Kuratowski-Mazurkiewicz theorem (on coverings of
simplices).

\medskip
With this goal in mind we use the theory of {\em quasitoric
manifolds} as introduced by  Davis and Januszkiewicz  in the seminal paper
\cite{Davis} and developed by many authors, see the monograph of
Buchstaber and Panov \cite{BuPan} (and the forthcoming,
considerably updated and revised new version \cite{newBuPan})).
Quasitoric manifolds offer more flexibility and versatility than
toric varieties, since they are easier to construct and their
geometric and algebraic topological properties are even more
closely related to combinatorics of simple polytopes.

\medskip
Another input came from the theory of projectives in simplicial
complexes and colorings of simple polytopes, as initiated by
Joswig in \cite{Joswig}. In particular we focus our attention to
the class of {\em $n$-colorable simple polytopes} (and some
generalizations) which appear to be particularly suitable as a
combinatorial framework for theorems of Lebesgue and KKM type.

\subsection{Coloring of simple polytopes}\label{sec:color-simple}

An $n$-dimensional convex polytope $P$ is \textit{simple} if the
number of codimension-one faces meeting at each vertex is $n$.
Codimension-one faces are called \textit{facets}. The following
inconspicuous lemma records one of the key properties of simple
polytopes.

\begin{lemma}
If $P$ is a simple polytope then two facets $F_1\neq F_2$ have a
non-empty intersection if and only if they share a common facet,
i.e.\ if $F_1\cap F_2$ is a face of $P$ of codimension $2$.
\end{lemma}\label{lemma:simple}
Suppose that $\{F_i\}_{i=1}^m$ is an enumeration of all facets of
a simple polytope $P^n$. A \textit{proper coloring} of $P^n$ by
$k$ colors is a map
 \begin{equation}\label{eqn:coloring}
  h: \{F_1, \dots, F_m\}\rightarrow [k]
 \end{equation}
(or a map $h : [m]\rightarrow [k]$) such that for each two
distinct facets $F_i\neq F_j$ if $F_i$ and $F_j$ are adjacent (in
the sense that they have a common facet) then $h (F_i)\neq h
(F_j)$.

In light of Lemma~\ref{lemma:simple} it is clear that $h$ is a
coloring of a simple polytope $P^n$  if and only if it is
a coloring of the graph on $[m]$ as the set of vertices, where
$(i,j)$ is an edge if and only if $F_i\cap F_j\neq\emptyset$. For
this reason the smallest number $k$ of colors needed for a proper
coloring of a simple polytope $P^n$ is called \textit{the
chromatic number} $\chi (P^n)$.

\medskip
It is immediate that $\chi (P^n)\geq n$ for any simple polytope
$P^n$. The chromatic number of a $2$-dimensional simple polytope
is clearly equal to $2$ or $3$, depending on the parity of the
number of its facets. By the {\em Four Color Theorem} we know that
the chromatic number of a $3$-dimensional polytope is either $3$
or $4$. However in general (for $n\geq 4$) it is far from being
true that $\chi (P^n)\leq n+1$. Actually one can easily produce
simple polytopes such that their chromatic number is exactly the
number of their facets. Examples include polytopes which arise as
polars of  cyclic polytopes $C^n (m)$, see \cite[Example~0.6,
p.11]{BuPan}.

\medskip
In spite of that the class of {\em $n$-colorable} $n$-dimensional
simple polytopes is quite large, with many interesting examples.
It is known that this class is closed for products
\cite[Construction~1.12, p.10]{BuPan} and connected sums
\cite[Construction~1.13, p.10]{BuPan}. From any given simple
polytope $P^n$ by truncation over all its faces we obtain a simple
polytope $Q^n$ such that $\chi (Q^n)=n$. The complete description
of this class is given by M.~Joswig in \cite{Joswig}, who proved
that a simple $n$-polytope $P^n$ admits an $n$-coloring if and
only if every $2$-face has an even number of edges. For this
reason an {\em $n$-colorable polytope} is sometimes referred to as
\textit{Joswig polytope}.

\begin{definition}\label{def:color-class}
Suppose that $P^n$ is an $n$-colorable simple polytope and let $h$
be an associated coloring function {\em (\ref{eqn:coloring})}. For
$0\leq k\le n$ let $I = \{i_1, i_2,\ldots, i_{n-k}\}\subset [n]$
be a collection of $(n-k)$ colors. We say that a $k$-dimensional
face $K$ of $P^n$ is in the $I$-color class if $I = \{h(F)\mid
K\subset F\}$.
\end{definition}

\begin{example}\label{exam:cube}{\em
The $n$-dimensional cube $I^n\subset \mathbb{R}^n$ is an
$n$-colorable simple polytope with colors corresponding to axes of
symmetry of pairs of opposite facets (coordinate axes). Similarly,
$I$-color classes of $k$-faces correspond to $(n-k)$-dimensional
coordinate subspaces of the ambient space $\mathbb{R}^n$. }
\end{example}

\begin{figure}[hbt]
\centering
\includegraphics[scale=0.3]{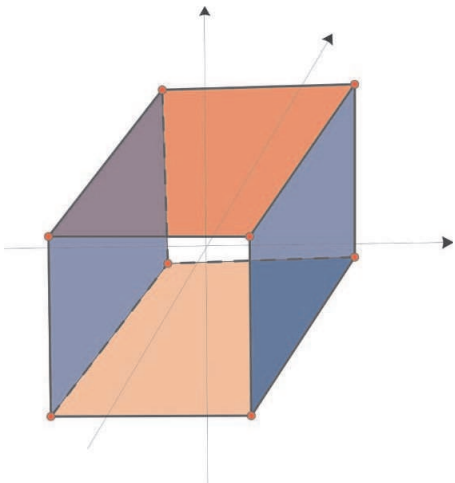}
\caption{Coloring of the cube $I^n$ (the case $n=3$).}
\label{fig:kocka}
\end{figure}

\subsection{Preliminaries on quasitoric manifolds}

A quasitoric manifold (originally a "toric manifold") is a
topological counterpart to the nonsingular projective toric
variety (of algebraic geometry). A smooth $2n$-dimensional
manifold $M^{2n}$ is a quasitoric manifold if it admits a smooth,
locally standard action of an $n$-dimensional topological torus
$T^n = (S^1)^n$, with an $n$-dimensional simple convex polytope
$P^n$ as the orbit space. Quasitoric manifolds were introduced by
Davis and Januszkiewicz in \cite{Davis} and developed by many
authors, see the monograph of Buchstaber and Panov \cite{BuPan}
and the forthcoming new version \cite{newBuPan} summarizing the
development of the theory in the last two decades.

\medskip
The facets $F_j$ of the polytope $P^n$ correspond to
$T^{n-1}$-orbits and the associated stabilizer groups define the
{\em characteristic map} (characteristic matrix) $\lambda: F_j
\mapsto T(F_j)$, where $T(F_j) = (\lambda_{ij})_{i=1}^n \in
\mathbb{Z}^n$ is a unimodular vector.

Conversely, each $n\times m$ {\em characteristic matrix} $\lambda
= (\lambda_{ij})$ produces a $2n$-dimensional quasitoric variety
$M^{2n}$ over a simple $n$-dimensional polytope $P^n$, provided
the column vectors $\lambda_j = (\lambda_{ij})$ satisfy the
condition that $\lambda_{j_1},\ldots,\lambda_{j_n}$ is a
$\mathbb{Z}^n$-basis for each choice of facets $F_{j_1},\ldots,
F_{j_n}$ having a common vertex.

\medskip

Following Davis-Januszkiewicz \cite[Theorem~4.14.]{Davis} there is
an isomorphism,
\begin{equation}\label{eqn:DJ-ideal-1}
   H^\ast(M^{2n}; \mathbb{Z}) \cong \mathbb{Z}[v_1,\ldots, v_m]/\langle I
   + J\rangle
\end{equation}
where $I$ is the {\em Stanley-Reisner} ideal of $P$ (generated by
monomials $v_{i_1}\ldots v_{i_k}$ such that $F_{i_1}\cap\ldots\cap
F_{i_k}=\emptyset$) and $J$ is the ideal generated by linear forms
which corresponds to the rows of the characteristic matrix
$(\lambda_{ij})$.

\subsection{Lyusternik-Schnirelmann method}

\begin{definition}\label{def:inessential}
For a given cohomology class $\omega\in H^\ast(X)$ we say that a
closed (open) subset $F\subset X$ is $\omega$-inessential (or
simply inessential if $\omega$ is clear from the context) if
$\omega$ is mapped to zero by the restriction map,
\[
H^\ast(X) \longrightarrow H^\ast(F).
\]
\end{definition}
\noindent The following well known `lemma' captures the essence of
the Lyusternik-Schnirelmann method.

\begin{lemma}\label{lemma:LS}
Assume that $\{X_i\}_{i=1}^n$ is a collection of closed (open)
subsets of a space $X$ and let $\{\omega_i\}_{i=1}^n$ be a
collection of cohomology classes in $H^\ast(X)$. If $X_i$ is
$\omega_i$-inessential for each $i=1,\ldots, n$ then $Z =
\cup_{i=1}^n~X_i$ is $\omega$-inessential where $\omega =
\omega_1\ldots\omega_n$.
\end{lemma}

The following proposition \cite[Lemma~3.2.]{Karasev} is the key
result connecting the {\em covering multiplicity} of a finite
family $\{Y_i\}_{i\in I}$ of subspaces of $Y$  with the cup-length
of the ring $H^\ast(Y)$. Recall that the covering multiplicity of
$\{Y_i\}_{i\in I}$ is $\leq k$ if for each $y\in Y$ the
cardinality of the set $\{i\in I \mid y\in Y_i\}$ is at most $k$.

\begin{proposition}\label{prop:crucial}
Suppose that a finite family $\mathcal{U} = \{U_i\}_{i=1}^N$ of
open subsets in a paracompact space $Y$ has covering multiplicity
at most $m$. Assume that for some $\omega\in H^\ast(Y)$ each of
the sets $U_i$ is $\omega$-inessential. Then the union
$\cup_{i=1}^N~U_i$ is $\omega^m$-inessential.
\end{proposition}

\begin{corollary}\label{cor:crucial}
Assuming that the cohomology theory satisfies a suitable
continuity condition (as the Alexandrov-\v{C}ech theory) the
Proposition~\ref{prop:crucial} is valid for finite closed
coverings of multiplicity $\leq m$.
\end{corollary}

\section{Colorful Lebesgue theorem}

Suppose that $P^n$ is an $n$-colorable simple polytope
(Section~\ref{sec:color-simple}) with $m$ facets $F_1, \ldots,
F_m$ and the corresponding coloring function (\ref{eqn:coloring}).
Let $e_1,\ldots, e_n$ be the standard basis of the lattice
$\mathbb{Z}^n$.

\begin{definition}\label{def:canonical}
The coloring {\em (\ref{eqn:coloring})} gives rise to a canonical
characteristic function $\lambda$ where $\lambda(F_i) = e_{h(i)}$.
The quasitoric manifold arising from this construction is referred
to as the {\em canonical quasitoric manifold} of the pair $(P^n,
h)$ or simply as a {\em canonical quasitoric manifold} associated
to the $n$-colorable simple polytope $P^n$.
\end{definition}

Suppose that $M^{2n}$ is the canonical quasitoric manifold
associated to an $n$-colorable simple polytope $P^n$. Let $\pi :
M^{2n}\rightarrow P^n$ be the corresponding projection map. For
each facet $F_i$ the set $M_i:=\pi^{-1}(F_i)$ is a codimension 2
submanifold of $M^{2n}$. Let $v_i\in H^2(M^{2n}; \mathbb{Z})$ be
the Poincar\'{e} dual of the fundamental class $[M_i]\in
H_2(M^{2n}; \mathbb{Z})$ (relative to some (omni)orientation on
$M^{2n}$).

According to the Davis-Januszkiewicz description of the
cohomological ring of $M^{2n}$ \cite[Theorem~4.14.]{Davis} there
is an isomorphism,
\begin{equation}\label{eqn:DJ-ideal}
   H^\ast(M^{2n}; \mathbb{Z}) \cong \mathbb{Z}[v_1,\ldots, v_m]/\langle I
   + J\rangle
\end{equation}
where $I$ is the {\em Stanley-Reisner} ideal of $P$ (generated by
monomials $v_{i_1}\ldots v_{i_k}$ such that $F_{i_1}\cap\ldots\cap
F_{i_k}=\emptyset$) and $J$ is the ideal generated by linear forms
$L_i$ where,
\begin{equation}\label{eqn:linear}
L_i(v_1,\ldots, v_m) = \sum_{h(j)=i}~v_j.
\end{equation}

\medskip The following proposition records some of the properties of
the cohomology ring of the canonical quasitoric manifold $M^{2n}$
associated to an $n$-colorable simple polytope $P^n$
(Definition~\ref{def:canonical}).

\begin{proposition}\label{prop:1234} \hfill
\begin{enumerate}
 \item[{\em (1)}] The product $v_iv_j$ of two distinct classes of the
same color is zero in $H^*(M^{2n}; \mathbb{Z})$.
 \item[{\em (2)}] The sum of all `classes of the same color'
 vanishes, $\sum_{h(i)=k}~v_i = 0$.
 \item[{\em (3)}] The square $v_i^2$ of any generator $v_i$ is
zero in $H^*(M^{2n}; \mathbb{Z})$.
  \item[{\em (4)}] Suppose that $\{F_{i_k}\}_{k=1}^n$ are all
  facets which share a common vertex $V$ of $P^{2n}$. Then,
  $(v_{i_1}+\ldots + v_{i_k})^n$ is a non-zero class in $H^{2n}(M^{2n}; \mathbb{Z})$.
\end{enumerate}
\end{proposition}

\textit{Proof.} The first observation is a direct consequence of
the fact that $v_iv_j\in I$ if $h(F_i)=h(F_j)$ and $i\neq j$. The
second property is just a restatement of the equation
(\ref{eqn:linear}) describing the ideal $J$. The property (3)
follows on multiplying the both sides of the equation $L_i = 0$ by
$v_i$. Finally (4) follows from the observation that
$v_{i_1}\ldots v_{i_n}$ is the fundamental cohomology class in
$H^{2n}(M^{2n}; \mathbb{Z})$ and the equality,
\[
(v_{i_1}+\ldots + v_{i_n})^n = n!\, v_{i_1}\ldots v_{i_n}.
\]
\hfill $\square$

The following result extends the Lebesgue theorem
(Theorem~\ref{thm:Lebesgue}) to the class of $n$-colorable simple
polytopes.

\begin{theorem}{\em (Colorful Lebesgue theorem)} \label{thm:color-Lebesgue}
Suppose that an $n$-colorable simple polytope $P^n$ is covered by
a family of closed sets $P^n = \cup_{i=1}^N~ X_i$ such that each
point $x\in P^n$ is covered by no more than $n$ of the sets $X_j$.
Then for some $i$, a connected component of $X_i$ intersects at
least two distinct facets of $P^n$ of the same color.
\end{theorem}

 \textit{Proof.} Without loss of generality we may assume that all
sets $X_i$ are connected. Indeed, the connected components of all
sets $X_j$ define a covering of $P^d$ which also satisfies the
conditions of the theorem. Let $M^{2n}$ be the canonical
quasitoric manifold over $P^n$ (Definition~\ref{def:canonical})
and let $\pi : M^{2n}\rightarrow P^n$ be the associated projection
map.

Assume (for contradiction) that each of the sets $X_j$ intersects
at most one facet of each of the colors $i\in [n]$. Given a vertex
$V$ of $P$, let $\{F_{i_k}\}_{k=1}^n$ be the collection of all
facets of $P^n$ incident to $V$ where  $h(F_{i_k})= h(i_k)=k$ for
the chosen coloring function (\ref{eqn:coloring}). By assumption
for each $k$ either $F_{i_k}\cap X_j =\emptyset$ (and
$\pi^{-1}(X_j)$ is automatically $v_{i_k}$-inessential) or
$F_i\cap X_j=\emptyset$ for each $F_i\neq F_{i_k}$ in the chosen
color class ($h(i)=k$). In the latter case $\pi^{-1}(X_j)$ is
$v_i$-inessential for each $i$ such that $h(i)=k$ and $F_i\neq
F_{i_k}$. Since the sum of all classes of the same color vanishes
(Proposition~\ref{prop:1234}) we conclude that $\pi^{-1}(X_j)$ is
$v_{i_k}$-inessential in this case as well.

Summarizing, we observe that $\pi^{-1}(X_j)$ is
$\omega$-inessential for each $j$ where $\omega = v_{i_1}+\ldots +
v_{i_n}$. It follows from Proposition~\ref{prop:crucial}
(Corollary~\ref{cor:crucial}) that $M^{2n} =
\cup_{j=1}^N~\pi^{-1}(X_j)$ is $\omega^n$-inessential which is in
contradiction with Proposition~\ref{prop:1234}.
 \hfill $\square$

\medskip

Theorem~\ref{thm:color-Lebesgue} extends the Lebesgue theorem
(Theorem~\ref{thm:Lebesgue}) to the class of all $n$-colorable
simple polytopes. Informally it says that if a collection
$\{X_i\}_{i=1}^N$ of closed subsets of $P^n$ has ``small
multiplicity'' (multiplicity $\leq n$) and sets of ``small
diameter'' ($X_i\cap F_j\neq\emptyset$ for at most one index $j$)
then it cannot be a covering of $P^n$.

\medskip
Karasev proved \cite[Theorem~4.2.]{Karasev} a very interesting
extension of Theorem~\ref{thm:Lebesgue} where he was able to show
that in a very precise sense the smaller is the multiplicity of
$\{X_i\}_{i=1}^N$ the bigger are the connected components of
$P^n\setminus \cup_{i=1}^N~X_i$. He obtained this result by
applying his method to the toric variety $(\mathbb{C}P^1)^n$ over
the cube $I^n$. Our objective is to extend this result to the
class of $n$-colorable simple polytopes.

\medskip
A `vertex class' $\omega\in H^2(M^{2n})$, associated to a vertex
$V\in P^n$, is by definition the sum $\omega = v_1+\ldots + v_n$
of all $2$-classes dual to facets $F_i$ incident to $V$.

\begin{theorem}\label{thm:quantitative-Lebesgue}
Suppose that $P^n$ is an $n$-colorable simple polytope, $M^{2n}$
its canonical quasitoric manifold, and $\pi : M^{2n}\rightarrow
P^n$ the associated projection map. Let $\omega = v_1+\ldots +
v_n$ be the $2$-dimensional `vertex class' associated to a vertex
$V\in P^n$. Suppose that ${\mathcal F} = \{X\}_{i=1}^N$ is a
finite family of closed subsets of $P^n$ such that each $X_i$
intersects at most one of the facets in each of the color classes.
If the covering multiplicity of ${\mathcal F}$ is at most $k\leq
n$ then there exists a connected component $Z$ of the set
$P^n\setminus\cup_{i=1}^N~X_i$ which is $\omega^{n-k}$-essential
in the sense that the restriction of the class $\omega^{n-k}$ on
$\pi^{-1}(Z)$ is non-trivial. Moreover, if $\mathcal{K}$ is the
collection of all $k$-faces $K$ of $P^n$ such that $Z\cap
K\neq\emptyset$ then $\mathcal{K}$ contains a collection of
$k$-faces of size at least  $2^{n-k}$ which are all in the same
$I$-color class for some $I=\{i_1,\ldots, i_{n-k}\}\subset [n]$
(Definition~\ref{def:color-class}).
\end{theorem}

\textit{Proof.}  For a chosen vertex $V\in P^n$ let
$\{F_i\}_{i=1}^n$ be the collection of all facets incident to $V$
(we assume that $h(F_i)=i$). If $v_i\in H^2(M^{2n}; \mathbb{Z})$
is the class associated to the facet $F_i$ then  (following
Proposition~\ref{prop:1234}) the class $\omega^n$ is non-zero
where $\omega = v_1+\ldots + v_n$. Moreover we observe that,
\begin{equation}\label{eqn:omega-k}
\omega^k = k! \sum_{J} v_J \neq 0
\end{equation}
where the sum is taken over all collections $J =\{j_1,\ldots,
j_{k}\}$ of colors of size $k$ and $v_J = v_{j_1}v_{j_2}\cdots
v_{j_{k}}$.

Simplifying the notation, from here on we say that $Y\subset P^n$
is $\omega^k$-inessential if the set $\pi^{-1}(Y)\subset M^{2n}$
is $\omega^k$-inessential. Assuming that  each $X_i$ intersects at
most one of the facets in each of the color classes we deduce (as
in the proof of Theorem~\ref{thm:color-Lebesgue}) that the set
$\cup_{i=1}^N~X_i$ and is $\omega^k$-inessential. Moreover
(assuming the cohomology is continuous) this holds also for some
small open neighborhood $U$ of $\cup_{i=1}^N~X_i$. It follows that
the restriction of $\omega^{n-k}$ on $\pi^{-1}(W)$ is non-trivial
where $W=P^n\setminus\cup_{i=1}^N~X_i$. Otherwise $W$ would be
$\omega^{n-k}$-inessential and $P^n = U\cup W$ would be
$\omega^n$-inessential (contradicting
Proposition~\ref{prop:1234}). Since $W$ is
$\omega^{n-k}$-essential the same holds for some connected
component $Z$ of $W$.

\medskip
In order to prove the second part of the theorem it will be
sufficient to show that there exists a collection of monomials
$v_J = v_{j_1}v_{j_2}\cdots v_{j_{n-k}}$ of size $\geq 2^{n-k}$
such that $Z$ is $v_J$-essential and all these monomials are in
the same $I$-color class (in the sense of
Definition~\ref{def:color-class}). Indeed, by the same argument as
before, if $Z$ is  $v_J$-essential then $Z\cap K\neq\emptyset$
where $K = F_{j_1}\cap\ldots\cap F_{j_{n-k}}$.

\medskip
In light of the fact that $Z$ is $\omega^{n-k}$-essential, by
inserting $n-k$ in the place of $k$ in the equality
(\ref{eqn:omega-k}) we observe that at least one of the monomials
$v_J = v_{j_1}v_{j_2}\cdots v_{j_{n-k}}$ must be non-zero in
$H^\ast(Z)$. Since (Proposition~\ref{prop:1234}),
\begin{equation}\label{eqn:replace}
v_{j_1} = \sum \{v_j \mid j\neq j_1 \mbox{ {\rm and} } h(j) =
h(j_1)\}
\end{equation}
we can replace the generator $v_{j_1}$ in the monomial $v_J$ by
the sum of the remaining generators $v_j$ of the same `color'
($h(j) = h(j_1)$). In other words we multiply both sides of the
equality (\ref{eqn:replace}) by the monomial $v_{j_2}\cdots
v_{j_{n-k}}$ and observe that on the right hand side there must
appear a monomial $v_J' = v_{j}v_{j_2}\cdots v_{j_{n-k}}$, in the
same $I$-color class as $v_J$, which is also non-zero in
$H^\ast(Z)$. This procedure can be continued for other indices
(generators) which guarantees that there exist at least $2^{n-k}$
different monomials in the same $I$-color class which are all
non-zero in $H^\ast(Z)$. \hfill $\square$

\begin{remark}{\em
The proof of Theorem~\ref{thm:quantitative-Lebesgue} shows that a
vertex $V$ of $P^n$ and the associated vertex class $\omega =
v_1+\ldots + v_n$ can be prescribed in advance. From here we
deduce that the vertex $V$ is certainly contained by one of the
$2^{n-k}$ $k$-dimensional faces in the same $I$-color class which
intersect $Z$. }
\end{remark}

\section{Colorful KKM-theorem}

In this section we prove a colorful version of
Knaster-Kuratowski-Mazurkiewicz `lemma' (Theorem~\ref{thm:KKM}).
The strategy is the same as in the previous section. We describe a
family of simple polytopes together with associated natural
quasitoric manifolds and show that each of them has a special
cohomology class $\omega$ such that $\omega^n\neq 0$.

\begin{definition}\label{def:special}
Suppose that a simple polytope $P^n$ can be colored by $(n+1)$
colors (in the sense of Section~\ref{sec:color-simple}) and for a
chosen coloring let $\{T_1, \ldots , T_k\} = h^{-1}(n+1)$ be the
collection of all facets colored by the color $n+1$. The polytope
$P^n$ is called {\em specially $(n+1)$-colorable} if the
associated coloring function $h: \{F_1, \dots, F_m\}\rightarrow
[n+1]$ has the property that all facets $\{T_i\}_{i=1}^k$ are
$n$-simplices.
\end{definition}

\begin{figure}[hbt]
\centering
\includegraphics[scale=0.20]{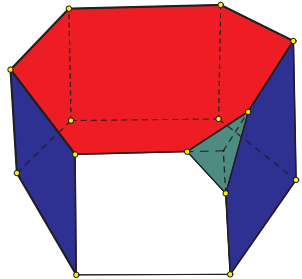}
\caption{A truncated $(n+1)$-colorable polytope (the case $n=3$).}
\label{fig:trunkacija}
\end{figure}

An immediate example of a special $(n+1)$-colorable polytope is
the standard simplex $\Delta^n$. A large class of such polytopes
is obtained by truncating $n$-colorable polytopes at an odd number
of (strongly separated) vertices (Figure~\ref{fig:trunkacija}).

\medskip
In the following definition we introduce a class of {\em canonical
quasitoric manifolds} associated to a $(n+1)$-colorable simple
polytope with a distinguished color (the color $n+1$).

\begin{definition}\label{def:sign}
Suppose that $P^n$ is a $(n+1)$-colorable, simple polytope. For
some enumeration $\{F_1, \dots, F_m\}$ of its facets let $h:
[m]\rightarrow [n+1]$ be a chosen coloring function. Let
$e_1,\ldots, e_n$ be the standard basis in $\mathbb{Z}^n$ and let
$e_\epsilon = \epsilon_1e_1+\ldots +\epsilon_ne_n$ where
$\epsilon_i\in \{-1,+1\}$. Define the characteristic function
$\lambda_\epsilon : \{F_1,\ldots, F_m\}\rightarrow \mathbb{Z}^n$
by the equation,
\begin{equation}\label{eqn:array}
\lambda_\epsilon(F_i) = \left\{
\begin{array}{ll}
e_{h(i)} & \mbox{ {\rm if} }\, i\neq n+1 \\
e_\epsilon  & \mbox{ {\rm if} }\, i = n+1
\end{array} \right.
\end{equation}
The quasitoric manifold $M^{2n}_\epsilon$ associated to the pair
$(P^n,\lambda_\epsilon)$ is called the {\em canonical quasitoric
manifold} of the $(n+1)$-colored polytope $P^n$ with the
distinguished color $n+1$ and the defining sign vector
$e_\epsilon$.
\end{definition}

By definition there are $2^n$ distinct canonical quasitoric
manifolds $M^{2n}_\epsilon$ associated to a $(n+1)$-colored
polytope $P^n$. From here on we select $e = -e_1-\ldots - e_n$ as
the preferred sign vector and denote the corresponding manifold by
$M^{2n}$.

\medskip
Proposition~\ref{prop:abcd} collects some of the properties of the
cohomology ring $H^\ast(M^{2n}; \mathbb{Z})$. This ring is, in
agreement with (\ref{eqn:DJ-ideal}), a quotient ring of
$\mathbb{Z}[v_1,\ldots,v_m]$ where $v_i$ is the $2$-dimensional
cohomology class associated to the facet $F_i$. For bookkeeping
purposes we modify (refine) this notation as recorded by the
following definition.

\begin{definition}\label{def:bookkeeping}
Following the notation of Definition~\ref{def:special}, let $t_i$
be the variable associated to the facet $T_i$. If
$\{F_{i\nu}\}_{\nu\in S_i}$ are facets colored by the color $i\in
[n]$ then the associated $2$-classes are $v_{i\nu}$. By assumption
$T_j$ is a simplex for each $j=1,\ldots, k$. It follows that for
each color $i\in [n]$ there is a unique facet $F_{ij}:=
F_{i\nu_j}$ adjacent to $T_j$. The associated dual cohomology
class is denoted by $v_{ij}$.
\end{definition}

\begin{proposition}\label{prop:abcd}  Suppose that $P^n$ is a specially
$(n+1)$-colorable polytope (Definition~\ref{def:special}) and let
$M^{2n}$ be the associated canonical quasitoric manifold
corresponding  the sign vector $e = -e_1-\ldots- e_n$
(Definition~\ref{def:sign}). Then the cohomology ring
$H^\ast(M^{2n}; \mathbb{Z})$, as described by {\em
(\ref{eqn:DJ-ideal})}, has the following properties.
\begin{enumerate}
\item[{\em (a)}] The Stanley-Reisner ideal $I$ of $P^n$ contains
all monomials $v_iv_j$ such that $h(i)=h(j)$ and $i\neq j$. In
particular $t_it_j\in I$ and $v_{i\nu_1}v_{i\nu_2}\in I$ for
$i\neq j$ and $\nu_1\neq\nu_2$.

\item[{\em (b)}] For each $i\in [n]$ there is a linear relation in
the cohomology ring $H^*(M^{2n}; \mathbb{Z})$,
 \begin{equation}\label{eqn:B}
t_1+\ldots + t_k = \sum_{\nu\in S_i}\, v_{i\nu}
 \end{equation}
\item[{\em (c)}] For each $j=1,\ldots, k$ there is a relation
(Definition~\ref{def:bookkeeping})
$$
t_j^2 =  t_jv_{ij}.
$$
\end{enumerate}
\end{proposition}

\textit{Proof.} The proof follows the same pattern as the proof of
Proposition~\ref{prop:1234}. For example on multiplying the
relation (\ref{eqn:B}) by $t_j$ one obtains $t_j^2 =
t_jv_{i\nu_j}$. Note that the equality (\ref{eqn:B}) is a
consequence of our choice of $e = -e_1-\ldots - e_n$ as the
preferred sign vector in the definition of the canonical
quasitoric manifold (Definition~\ref{def:sign}).
 \hfill $\square$

\begin{proposition}\label{prop:fundamental}
The class $t_j^n$ is equal to the fundamental cohomology class of
$M^{2n}$ for each $j=1,\ldots, k$.
\end{proposition}

\textit{Proof.} It follows from Proposition~\ref{prop:abcd} that,
\begin{equation}\label{eqn:t-na-n}
t_j^n =
t_j^2t_j^{n-2}=v_{1j}t_j^{n-1}=v_{1j}v_{2j}t_j^{n-2}=\ldots =
v_{1j}v_{2j}\ldots v_{n-1,j}t_j.
\end{equation}
If $\pi : M^{2n}\rightarrow P^n$ is the projection map than $t_j$
is the dual to the fundamental homology class of the (properly
oriented) manifold $\pi^{-1}({T_j})$, similarly $v_{ij}$ is dual
to the fundamental homology class of $\pi^{-1}(F_{ij})$. All these
manifolds intersect transversely and their intersection is a
single point. It follows from Poincar\'{e} duality that
$v_{1j}v_{2j}\ldots v_{n-1,j}t_j$ is the dual of a generator of
$H_0(M^{2n}; \mathbb{Z})$ so both $v_{1j}v_{2j}\ldots
v_{n-1,j}t_j$ and $t_j^n$ are fundamental classes of $M^{2n}$.
\hfill $\square$

\begin{theorem}\label{thm:color-KKM}
{\em (Colorful KKM theorem)} Let $P^n$ be a specially
$(n+1)$-colorable polytope in the sense of
Definition~\ref{def:special}. Suppose that $P^n$ is covered by a
family of closed sets $P^n = \cup_{i=1}^N~ X_i$ with the covering
multiplicity $\leq n$ (i.e.\ each point $x\in P^n$ is covered by
no more than $n$ of the sets $X_j$). Then there exists $i\in [N]$
and a connected component $Y_i$ of $X_i$ such that among the faces
of $P^{n}$ intersected by $Y_i$ are facets of all $n+1$ colors.
\end{theorem}

 \textit{Proof.} As in the proof of
Theorem~\ref{thm:color-Lebesgue} we are allowed to assume that all
sets $X_i$ are connected. Let $M^{2n}$ be the canonical quasitoric
manifold associated to $P^n$ corresponding to the sign vector $e =
-e_1-\ldots- e_n$ (Definition~\ref{def:sign} and
Proposition~\ref{prop:abcd}). Let $t = t_1+t_2+\ldots + t_k$ be
the sum of all $2$-classes corresponding to simplicial facets
$T_j$ of the polytope $P^n$.

\medskip
If $X_j\cap F_{i\nu}=\emptyset$ for each facet $F_{i\nu}$ of color
$\nu$ then $X_j$ is $v_{i\nu}$-inessential for each $\nu\in S_i$.
We deduce from Proposition~\ref{prop:abcd}~(b) that $X_j$ is
$t$-inessential as well and by Proposition~\ref{prop:crucial}
(Corollary~\ref{cor:crucial}) we know that $M^{2n} =
\cup_{j=1}^N~\pi^{-1}(X_j)$ is $t^n$-inessential. This contradicts
the fact that the class (Proposition~\ref{prop:fundamental})
$$
t^n = (t_1+\ldots + t_k)^n = t_1^n+\ldots + t_k^n = k t_1^n
$$
is a non-trivial element of $H^{2n}(M^{2n}; \mathbb{Z})$.
 \hfill $\square$

\medskip
It is certainly possible to refine Theorem~\ref{thm:color-KKM}
along the lines of Theorems~2.1. and 4.2. in \cite{Karasev} and
our Theorem~\ref{thm:quantitative-Lebesgue}. The following
corollary of the proof of Theorem~\ref{thm:color-KKM} summarizes
the cohomological content of such a result.

\begin{theorem}\label{thm:color-KKM-bis}
Let $P^n$ be a specially $(n+1)$-colorable polytope
(Definition~\ref{def:special}). Suppose that $P^n$ is covered by a
family of closed sets $P^n = \cup_{i=1}^N~ X_i$ with the covering
multiplicity $k \leq n$ and there is no $X_i$ intersecting some
$n+1$ distinct colored facets. Then there exists a connected
component $W$ of $P^n\setminus \cup_{i=1}^N~X_i$ which is
$t^{n-k}$-essential in the sense that the restriction of the class
$t^{n-k}$ on $\pi^{-1}(W)\subset M^{2n}$ is non-trivial. Moreover,
if $\mathcal{K}$ is the collection of all $k$-faces $K$ of $P^n$
such that $W\cap K\neq\emptyset$ then $\mathcal{K}$ contains a
 $k$-skeleton of some simplicial face $T_i$ and at least ${n \choose k}$
 $k$-faces of $P^n$ not contained in $T_i$.
\end{theorem}

 \textit{Proof.} The proof uses the same arguments as the first
 half of the proof of Theorem~\ref{thm:quantitative-Lebesgue} so we omit the
 details.  \hfill $\square$

\section{General Polytopes}

In this section we briefly address the case of general (not
necessarily simple) polytopes. We use the fact that after
truncations along all the faces of a polytope $P^n$ we obtain a
Joswig polytope $\overline{P}^n$. Indeed, facets $F_K$ of
$\overline{P}^n$ are naturally indexed by faces $K$ of $P^n$ and a
proper coloring of $\overline{P}^n$ by $n$ colors is defined by
$h(F_K) = \mbox{\rm dim}(K)$.

\begin{theorem}  Let a  polytope $P^n$ be covered by a
family of closed sets $\left\{X_i\right\}^N_{i=1}$ with covering
multiplicity at most $n$. Then some connected component of $X_i$
intersects at least two different $k$-faces of $P^n$ (for some
$k$).
\end{theorem}

\textit{Proof.} Let $\overline{P}^n$ be the total truncation of
$P_n$ such that $\partial \overline{P}^n$ lies in $\varepsilon$
neighborhood of $\partial P^n$, where $\varepsilon >0$ is a
sufficiently small, positive number. Observe that the restriction
of the family $\left\{X_i\right\}^N_{i=1}$ to $\overline{P}^n$ is
a covering of $\overline{P}^n$ by closed subsets. Theorem
\ref{thm:color-Lebesgue} implies that some connected component of
$X_i$ must intersect at least two facets $F_{K_1}$ and $F_{K_2}$
of $\overline{P}^n$, corresponding to $k$-faces $K_1$ and $K_2$ of
$Q^n$. The result follows by a limiting argument (when
$\varepsilon\rightarrow 0$). \hfill $\square$

\section{Colorful Hex Theorem}

As illustrated in previous sections quasitoric manifolds are a
very useful tool for analyzing various generalizations of the
Lebesgue and KKM theorem. Karasev observed in
\cite[Theorem~4.3.]{Karasev} that the Lyusternik-Schnirelmann
method is equally useful for the proof of the $n$-dimensional {\em
Hex theorem} \cite{Hex}.

\medskip
Suppose that $\{A_{i1}, A_{i2}\}_{i=1}^n$ is some labelling of
pairs of opposite facets of the $n$-dimensional cube $I^n$. The
Hex Theorem claims that for each covering $I^n = \cup_{i=1}^n~X_i$
of the $n$-cube $I^n$ by closed sets there exists an index $i$ and
a connected component $Y_i$ of $X_i$ such that both $Y_i\cap
A_{i1}\neq\emptyset$ and $Y_i\cap A_{i2}\neq\emptyset$.

\medskip
At first site this result is an immediate consequence and a very
special case of Theorem~\ref{thm:Lebesgue}, however on closer
inspection we see that this is not the case. Indeed in the Hex
Theorem we can prescribe in advance a matching between closed sets
$X_i$ and the corresponding pairs $\{A_{i1}, A_{i2}\}$ of opposite
facets of $I^n$.

\medskip
For this reason we formulate and prove a Hex analogue of
Theorem~\ref{thm:color-Lebesgue} which contains
\cite[Theorem~4.3.]{Karasev} as a special case.

\begin{theorem}{\em (Colorful Hex theorem)} \label{thm:color-Hex}
Suppose that $P^n$ is an $n$-colorable simple polytope and let $h
: [m]\rightarrow [n]$ be a selected coloring function which
associates to each facet $F_j$ the corresponding color $h(F_j) =
h(j)$. Let $V$ be a vertex of $P^n$ and let
$\{F_{\nu_i}\}_{i=1}^n$ be the collection of all facets of $P^n$
which contain $V$ such that $h(\nu_i)=i$.

Suppose that the polytope $P^n$ is covered by a family of $n$
closed sets $P^n = \cup_{i=1}^n~X_i$. Then for some $i$, a
connected component of $X_i$ intersects both $F_{\nu_i}$ and some
other facet $F_j$ of $P^n$ colored by the color $i$.
\end{theorem}

 \textit{Proof.} Let $M^{2n}$ be the canonical quasitoric
manifold associated to the simple $n$-colorable polytope $P^n$
with the chosen coloring $h : [m] \rightarrow [n]$ and let $\pi :
M^{2n} \rightarrow P^n$ be the corresponding projection map. Let
$v_j$ be the class associated to the facet $F_j$, in particular
$v_{\nu_i}$ is the class associated to the facet $F_{\nu_i}$.

 \medskip
It is sufficient to show that there exists a color $i$ such that
$X_i$ is $v_{\nu_i}$-essential in the sense that the restriction
of the class $v_{\nu_i}$ on $\pi^{-1}(X_i)$ is non-zero. Indeed,
in this case some connected component $Y_i$ of $X_i$ would be also
$v_{\nu_i}$-essential. From here it would follow that for some
$v_j\neq v_{\nu_i}$ such that $h(j)=i$ the set $Y_i$ would be
$v_i$-essential as well. Otherwise, in light of the relation,
$$
\sum_{h(j)=i} v_j =0
$$
we would deduce that $Y_i$ is NOT $v_{\nu_i}$-essential. This
observation would complete the proof since $Y_i$ would certainly,
under these condition, have a non-empty intersection with both
$F_{\nu_i}$ and $F_i$.

\medskip
Let us assume now that $X_i$ is $v_{\nu_i}$-inessential for each
$i=1,\ldots, n$. It follows from the Lyusternik-Schnirelmann
Lemma~\ref{lemma:LS} that $M^{2n} = \cup_{i=1}^n~\pi^{-1}(X_i)$ is
$\omega$-inessential where $\omega = v_{\nu_1}v_{\nu_2}\cdots
v_{\nu_n}$. This contradicts the fact that $\omega$ is the
fundamental cohomology class of $M^{2n}$.
 \hfill $\square$

\subsection{A generalized Game of Hex}

According to Wikipedia, \url{http://en.wikipedia.org/wiki/Hex_%28board_game%29}, see also
\cite{Hex, MZB}, it was {\em John Nash} who originally proved that
the game of Hex cannot end in a tie. Moreover, he is attributed to
be the first who observed that the first player has the winning
strategy for the game of Hex on the usual (rhombic) game board.

\medskip
Here, as an application of Theorem~\ref{thm:color-Hex}, we
describe a fairly general version of the Game of Hex \cite{Hex}
played by $n$ players which also cannot end in an undecided
position.

\begin{figure}[hbt]
\centering
\includegraphics[scale=0.80]{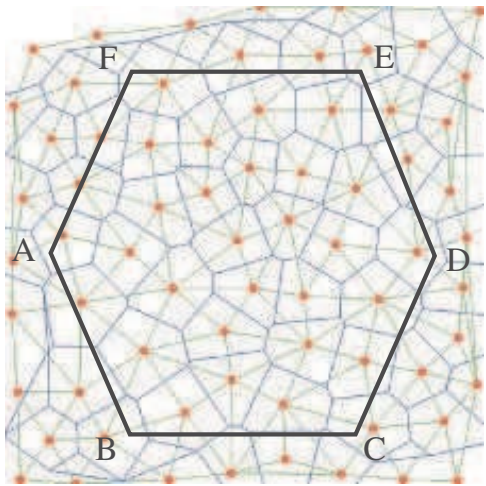}
\caption{A hexagonal Voronoi checkerboard.} \label{fig:Voronoi}
\end{figure}

Let $S\subset \mathbb{R}^n$ be a finite set of points and let
$\{V_x\}_{x\in S}$ the associated Voronoi partition of
$\mathbb{R}^n$. Let $P^n$ be a simple $n$-colorable polytope with
$m$ facets and an associated coloring function $h : [m]
\rightarrow [n]$. Choose a vertex $V$ of $P^n$ and let
$\{F_{\nu_i}\}_{i=1}^n$ be the collection of all facets containing
the vertex $V$ such that $h(\nu_i)=i$.

\medskip
There are $n$ players $J_1, \ldots , J_n$ (each player $J_i$ is
assigned the corresponding color $i$). The first player chooses a
point $x_1\in S$ and colors the corresponding Voronoi cell
$V_{x_1}$ by the color $1$. The second player chooses a point
$x_2\in S\setminus\{x_1\}$ and colors the Voronoi cell $V_{x_2}$
by color $2$, etc. After the first round of the game the first
player chooses a point $x_{n+1}\in S\setminus\{x_1,\ldots, x_n\}$,
etc.

\medskip
The game continues until one of the players (say the player $J_i$)
creates a connected monochromatic set of cells (all of color $i$)
which connect the facet $F_{\nu_i}$ with one of the facets $F_j$
such that $h(j)=i$.  Alternatively the game ends if there are no
more points in $S$ to distribute among players.

\medskip
An easy application of Theorem~\ref{thm:color-Hex} shows that the
game will always be decided i.e.\ sooner or later one of the
players will win the game.

\medskip
Figure~\ref{fig:Voronoi} illustrates the simplest new case of the
game. The red player chooses edges $AB, CD, EF$ of the hexagon,
while remaining edges belong to the blue player. The red player
tries to connect the edge $AB$ with either the edge $CD$ or $EF$,
similarly the blue player tries to connect the edge $AF$ with
either $CD$ or $EF$. As predicted by the Colorful Hex Theorem
sooner or later one of the players will achieve her goal.


\end{document}